\documentclass[11pt]{amsart} 
\usepackage{amssymb,epsfig}

\newtheorem{thm}{Theorem}
\newtheorem{lem}[thm]{Lemma}

\newtheorem{cor}[thm]{Corollary}
\newtheorem{conj}[thm]{Conjecture}

\theoremstyle{definition}

\newenvironment{pf}{\par\noindent \emph{Proof}.\ }{\medskip}

\title{On polynomial digraphs\footnote{
Work partially supported by the Ministerio de Ciencia y
Tecnolog\'ia under projects BFM2003-00368 and MTM2004-01728 and
Ministerio de Ciencia y Tecnolog\'ia and by the Generalitat de
Catalunya under project 2005 SGR 00692}}
\author{Josep M. Brunat \& Antonio Montes}
\address{Departament de Matem\`atica Aplicada II,
Universitat Polit\`ecnica de Catalunya, Spain.
 \{josep.m.brunat, antonio.montes\}@upc.edu,
http://www-ma2.upc.edu/$\sim$montes.}

\date{Mars, 2004.}
\begin{document}
%
\setcounter{page}{1}

\maketitle

\begin{abstract}
  Let $\Phi(x,y)$ be a bivariate polynomial with complex coefficients.
  The zeroes of $\Phi(x,y)$ are given a combinatorial structure by
  considering them as arcs of a directed graph $G(\Phi)$. This paper
  studies some relationship between
  the polynomial $\Phi(x,y)$ and the structure of $G(\Phi)$.

 \end{abstract}


\section{Introduction}

Let $\Phi(x,y)\in\mathbb{C}[x,y]$ be a bivariate polynomial with
complex coefficients and let $I$ be the ideal generated by
$\Phi(x,y)$.  The variety $\mathbb{V}(I)$ of $I$ is the set of ordered
pairs $(u,v)\in\mathbb{C}^2$ such that $\Phi(u,v)=0$. We give a
combinatorial structure to $\mathbb{V}(I)$ by taking its elements as
arcs of a digraph $G(\Phi)$, and we explore the relationship between
the polynomial $\Phi(x,y)$ and the digraph $G(\Phi)$.

The digraphs $G(\Phi)$ were introduced in~\cite{BrLa} for symmetric
polynomials $\Phi(x,y)$ under the name of \emph{Galois graphs}.
Lengths of walks, distances and cycles were described in terms of
$\Phi(x,y)$. Also, when the coefficients of $\Phi(x,y)$ belong to a
field $k$ and $\overline{k}$ is the algebraic closure of $k$, the
action of the Galois group $G(\overline{k}/k)$ on $G(\Phi)$ was
studied (and this was the motivation for the name of Galois graphs).

Here we give a general overview of the topic for non necessarily
symmetric polynomials and some concrete families of polynomials
are studied.  In~\cite{BrMo2} we adopt an algebraic approach to the
problem of deciding which polynomials produce a given graph as a
connected component of $G(\Phi)$.

Let us fix some notation. In this paper digraphs are allowed to be
infinite, and to have multiple arcs and loops. A \emph{graph} is a
digraph without loops nor multiple arcs such that for each arc $(u,v)$
there exists the arc $(v,u)$. The arcs $(u,v),(v,u)$ form the
\emph{edge} $uv$ of the graph.  Let $u$ be a vertex of a digraph $D$.
The \emph{strong (connected) component} of $u$ is the subdigraph of
$D$ induced by $u$ and the set of vertices $v$ such that there exist a
directed path from $u$ to $v$ and a directed path from $v$ to $u$.
The \emph{underlying graph} of a digraph $D$ is the graph obtained by
taking as edges the set of $uv$ with $(u,v)$ an arc of $D$ and $u\ne
v$. The \emph{(weakly connected) component} of a vertex $u$ in a
digraph $D$ is the subdigraph of $D$ induced by $u$ and all vertices
of the strong component of $u$ in the underlying graph of $D$. Note
that, in a graph, components and strong components coincide.  For
undefined concepts about graph theory, we refer
to~\cite{ChLe,Watkins2}. The component of $u$ in $G(\Phi)$ is
denoted by $G(\Phi,u)$ and the strong component by
$\vec{G}(\Phi,u)$.

In a monomial $cx^iy^j$, $c\ne 0$, the non negative integer $i$ is
called the \emph{partial degree} respect to $x$ (analogously for $j$
and $y$).  The \emph{total degree} or \emph{degree} of the monomial is
the integer $i+j$. The \emph{partial degree} respect to $x$ of a
bivariate polynomial $\Phi(x,y)$ is the maximum of the partial degrees
of its monomials; and the \emph{total degree} or \emph{degree} of
$\Phi(x,y)$ is the maximum of the degrees of its monomials. If
$\Phi(x,y)=\Phi_1(x,y)^{n_1}\cdots \Phi_k(x,y)^{n_k}$ is the expression of
$\Phi(x,y)$ as a product of irreducible polynomials over $\mathbb{C}$, the
\emph{radical} of $\Phi(x,y)$ is the polynomial
$\mathrm{rad\,}\Phi(x,y)=\Phi_1(x,y)\cdots \Phi_k(x,y)$. A polynomial is
\emph{radical} if $\mathrm{rad\,}\Phi(x,y)=\Phi(x,y)$.  We refer to
\cite{CoLiSh,CoLiSh2} for any other undefined concept about polynomials.

The paper is organized as follows. In the rest of this section we
give an informal description of its contents. In the next section
we define the digraph $G(\Phi)$. We shall see that some natural
conditions on $\Phi(x,y)$ (such as $\Phi(x,y)$ to be radical,
$\Phi(u,y)\ne 0$ for all $u\in\mathbb{C}$, etc.)  can be assumed.
Polynomials satisfying such conditions are called \emph{standard
polynomials}.  If $\Phi(x,y)$ is a standard polynomial, then
$G(\Phi)$ has only a finite number of loops and multiple arcs.
Moreover, all vertices have finite indegree and all of them, but a
finite number, have the same indegree; analogously for outdegrees.
The components (resp. strong components) of $G(\Phi)$ containing
these vertices, multiple arcs and loops are called \emph{singular
components} (resp. \emph{singular strong components}).

In Section~3 we show that every finite strongly connected $d$-regular
digraph is isomorphic to a strong component of $G(\Phi)$ for an appropriate
$\Phi(x,y)$.  Nevertheless, the construction produces an infinite
component. A digraph is called
\emph{polynomial} if it is isomorphic to $G(\Phi)$ or to a non
singular component of $G(\Phi)$.

In Section~4 we will see that Cayley digraphs on the additive and
multiplicative groups of $\mathbb{C}$ are polynomial and we give the
corresponding polynomial $\Phi(x,y)$. This implies that directed and
undirected cycles, finite complete graphs $K_d$, finite bipartite complete
graphs $K_{d,d}$ and, in general, circulant digraphs are polynomial.

In Section~5 we study polynomials of partial degree one in each
indeterminate. It is shown that all non-singular components are
isomorphic, and a characterization of polynomials of partial degree
one which give directed $n$-cycles as non-singular components is
given, as well as those giving infinite directed paths. By relating these
polynomials to the group of linear fractional transformations, we
prove that Cayley digraphs on dihedral groups and on the groups of
symmetries of regular polyhedra are polynomial.

In Section~6 we consider symmetric polynomials of total degree 2. In
this case the structure of $G(\Phi)$ is also completely determined.
In particular, all non-singular components are isomorphic.

Finally, the results obtained here and the discussions in~\cite{BrMo2}
give support to the conjecture stated in the last section.

\section{The digraph of a polynomial}
Let
$$
\Phi(x,y)=\sum_{i=0}^d a_i(x)y^i=\sum_{j=0}^e b_j(y)x^j
$$
be a polynomial with complex coefficients. The digraph $G(\Phi)$
has $\mathbb{C}$ as set of vertices and an ordered pair $(u,v)$ is an
arc of multiplicity $m$ if $v$ is a root of multiplicity $m$ of the
polynomial $\Phi(u,y)$. If $\Phi(u,y)$ is the zero polynomial for some
$u\in\mathbb{C}$ the multiplicity of all arcs $(u,v)$,
$v\in\mathbb{C}$, is taken to be 1, and the vertex $u$ is called a
\emph{source universal vertex}. Note that the source universal vertices
are the roots of the polynomial
$A(x)=\mathrm{gcd}(a_0(x),\ldots,a_d(x))$. Analogously, a vertex $v$
is called a \emph{sink universal vertex} if $\Phi(x,v)$ is the zero
polynomial; the sink universal vertices are the roots of the
$B(y)=\mathrm{gcd}(b_0(y),\ldots,b_e(y))$. In the following we assume
that $G(\Phi)$ has no universal vertices, or, equivalently that the
polynomials $A(x)$ and $B(y)$ are constant. If $\Phi(x,y)$ is constant
then $G(\Phi)$ is the complete or the null digraph on $\mathbb{C}$
depending on whether the constant is zero or not. Therefore, we can
also assume from now on that $e,d\ge 1$ and that $a_d(x)$ and $b_e(y)$
are non zero polynomials.

The structure of $G(\Phi)$ is given by the structure of its components.
To study the components of $G(\Phi)$, it is useful to put aside
some special cases.

Consider the discriminant
$$
D(x)=\mathrm{Resultant}(\Phi(x,y), \, \Phi'_y(x,y), \, y) \,,
$$
where $\Phi'_y(x,y)$ denotes the partial derivative of $\Phi(x,y)$
with respect to $y$.  We have that $D(x)=0$ if and only if $\Phi(x,y)$
has a multiple factor of positive degree in $y$, say
$\Phi(x,y)=\Phi_1(x,y)^k\Phi_2(x,y)$ with $k\ge 2$.  In this case,
the digraphs $G(\Phi)$ and $G(\mathrm{rad\,}\Phi)$ differ only in the
multiplicity of the arcs corresponding to the factor $\Phi_1(x,y)$.
Therefore, we can assume that $D(x)$ is not the zero polynomial.
Analogously, we can assume that
$$
E(y)=\mathrm{Resultant}(\Phi(x,y), \, \Phi'_x(x,y), \, x) \,,
$$
is not the zero polynomial.

Note that all vertices have outdegree at most $d$ and a vertex $u$ has
outdegree $<d$ if and only if $a_d(u)=0$. Analogously, all vertices
have indegree at most $e$ and a vertex $u$ has indegree $<e$ if and
only if $b_e(u)=0$. The roots of $a_d(x)$ are called
\emph{out-defective vertices} and the roots of $b_e(y)$ are called
\emph{in-defective vertices}.

Let $D(x)\ne 0$ and $A(x)=1$.  The leading coefficients of $\Phi(x,y)$
and $\Phi'_y(x,y)$ as polynomials in the indeterminate $y$ are
$a_d(x)$ and $d a_d(x)$ respectively, so $a_d(x)$ is a factor of
$D(x)$.  Thus, the out-defective vertices are roots of $D(x)$.  If $u$
is the origin of a multiple arc, then $\Phi(u,y)$ has a multiple root.
Therefore, $D(u)=0$. Conversely, if $D(u)=0$, then either $a_d(u)=0$
or $\Phi(u,y)$ has a multiple root, i.e. $u$ is an out-defective
vertex or it is the origin of a multiple arc. We conclude that the roots
of $D(x)$ are the out-defective vertices and the origins of multiple
arcs. Analogously, the roots of $E(y)$
are the in-defective vertices and the
ends of multiple arcs.

The vertices with a loop are the roots of $L(x)=\Phi(x,x)$. There is a
loop at each vertex if and only if $L(x)$ is the zero polynomial,
which means that $\Phi(x,y)$ admits a factorization
$\Phi(x,y)=(y-x)^k\Phi_1(x,y)$ with $k\ge 1$ and $\Phi_1(x,y)$ not
divisible by $y-x$. The structure of $G(\Phi)$ is then completely
determined by the structure of $G(\Phi_1)$, and thus it is not a
restriction to assume that $L(x)$ is not the zero polynomial.

A polynomial $\Phi(x,y)$ is \emph{standard} if it is non constant
($G(\Phi)$ is neither the complete or the null graph), $A(x)B(y)$ is
constant (there are not universal vertices), $D(x),E(x)\ne 0$ (the
polynomial $\Phi(x,y)$ is a radical polynomial), and $L(x)$ is not the
zero polynomial ($G(\Phi)$ has no loops at every vertex).  For a
standard polynomial $\Phi(x,y)$, the roots of $S(x)=L(x)D(x)E(x)$ are
called \emph{singular vertices}. They are the vertices with a loop,
vertices which are origin or end of multiple arcs, and defective
vertices. A \emph{singular component} (resp. \emph{singular strong
  component}) is a component (resp. strong component) which contains
some singular vertex.  Note that only a finite number of singular
components (strong components) exists.  We denote by $G(\Phi)^*$ the
digraph obtained from $G(\Phi)$ by removing all its singular components.

Note that different polynomials $\Phi(x,y)$ can give
isomorphic digraphs $G(\Phi)$, as stated in the following lemma of
straightforward proof.

\begin{lem}
\label{iso}
Let $\Phi(x,y)\in\mathbb{C}[x,y]$ and $a,b,c\in\mathbb{C}$ with
$a,c\ne 0$. If $\Psi(x,y)=c\Phi(ax+b,ay+b)$, then the mapping
$u\mapsto au+b$ is an isomorphism from $G(\Psi)$ to $G(\Phi)$.
\end{lem}

\section{Finite strong components}

Our immediate goal is to show that every finite strongly connected
$d$-regular digraph can be seen as a strong component of $G(\Phi)$ for
some appropriate $\Phi(x,y)$. We need the following Lemma:
\begin{lem}
\label{factorization}
Every $d$-regular digraph admits a 1-factorization.
\end{lem}
Lemma~\ref{factorization} can be proved by using that every
regular graph of even degree admits a
2-factorization~\cite{Petersen}. See also~\cite{Fabrega} for a detailed proof.

\begin{thm}
\label{strong}
Let $D=(V,E)$ be a finite strongly connected $d$-regular  digraph
of order $n\ge 2$ with $V\subset\mathbb{C}$.
Then there exists a polynomial $\Phi(x,y)$
such that $D$ is a strong component of $G(\Phi)$.
\end{thm}
\begin{pf}
  Lemma~\ref{factorization} ensures that $D$ admits a
  $1$-factorization. Let $F_1,\ldots,F_d$ be the set of arcs of the
  1-factors.  For each $i\in[d]=\{1,\ldots,d\}$ let
  $L_i(x)$ be the interpolation polynomial such that $L_i(u)=v$ for
  each $(u,v)\in F_i$. In this way, the vertices adjacent from $u\in
  V$ are $L_1(u),\ldots,L_d(u)$. Define $\Phi(x,y)=(y-L_1(x))\cdots
  (y-L_d(x))$. In $G(\Phi)$, the vertex $u$ is also adjacent to
  $L_1(u),\ldots, L_d(u)$. Therefore $D=\vec{G}(\Phi,u)$. \qed
\end{pf}

\begin{figure}[t]
\begin{center}
\begin{picture}(0,0)%
\epsfig{file=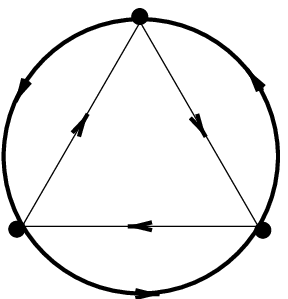}%
\end{picture}%
\setlength{\unitlength}{2486sp}%
\begingroup\makeatletter\ifx\SetFigFont\undefined%
\gdef\SetFigFont#1#2#3#4#5{%
  \reset@font\fontsize{#1}{#2pt}%
  \fontfamily{#3}\fontseries{#4}\fontshape{#5}%
  \selectfont}%
\fi\endgroup%
\begin{picture}(2136,2455)(1408,-2411)
\put(2416,-121){\makebox(0,0)[lb]{\smash{\SetFigFont{8}{9.6}{\familydefault}{\mddefault}{\updefault}$1$}}}
\put(1471,-2138){\makebox(0,0)[lb]{\smash{\SetFigFont{8}{9.6}{\familydefault}{\mddefault}{\updefault}$2$}}}
\put(3368,-2145){\makebox(0,0)[lb]{\smash{\SetFigFont{8}{9.6}{\familydefault}{\mddefault}{\updefault}$3$}}}
\end{picture}
\caption{A factorization of the digraph $K_3$}
\label{figk3}
\end{center}
\end{figure}

In the proof of Theorem~\ref{strong}, the polynomials $L_i(x)$
can have degree $n-1$, so $\Phi(x,y)$ can have degree $d$ in $y$ and degree
$d(n-1)$ in $x$. If the given digraph $D$ has order $n\ge 3$, then
$d(n-1)>d$ and the component $G(\Phi,u)$ is infinite, while the
strong component $\vec{G}(\Phi,u)=D$ is finite.  For instance, take
$D=K_3$, the complete symmetric digraph of order 3,
and choose $1,2$ and $3$ as
the vertices of $D$. The digraph $D$ is 2-regular and admits the
factorization $F_1,F_2$ where $F_1=\{(1,2),(2,3),(3,1)\}$ and
$F_2=\{(1,3),(3,2),(2,1)\}$. Figure~\ref{figk3} shows the
factorization of $K_3$; the arcs of $F_1$ are the thick ones.

\begin{figure}[h]
\begin{center}
\begin{picture}(0,0)%
\epsfig{file=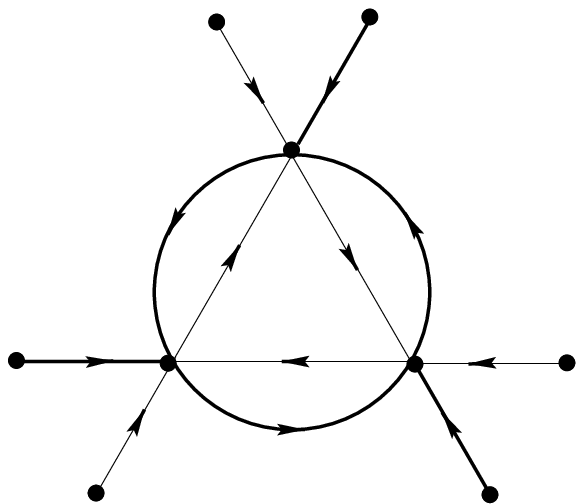}%
\end{picture}%
\setlength{\unitlength}{2030sp}%
\begingroup\makeatletter\ifx\SetFigFont\undefined%
\gdef\SetFigFont#1#2#3#4#5{%
  \reset@font\fontsize{#1}{#2pt}%
  \fontfamily{#3}\fontseries{#4}\fontshape{#5}%
  \selectfont}%
\fi\endgroup%
\begin{picture}(5397,4663)(369,-4079)
\put(3061,-631){\makebox(0,0)[lb]{\smash{\SetFigFont{6}{7.2}{\familydefault}{\mddefault}{\updefault}$1$}}}
\put(4486,-2978){\makebox(0,0)[lb]{\smash{\SetFigFont{6}{7.2}{\familydefault}{\mddefault}{\updefault}$3$}}}
\put(1711,-2978){\makebox(0,0)[lb]{\smash{\SetFigFont{6}{7.2}{\familydefault}{\mddefault}{\updefault}$2$}}}
\put(4013,419){\makebox(0,0)[lb]{\smash{\SetFigFont{6}{7.2}{\familydefault}{\mddefault}{\updefault}$7/3$}}}
\put(1913,419){\makebox(0,0)[lb]{\smash{\SetFigFont{6}{7.2}{\familydefault}{\mddefault}{\updefault}$2/3$}}}
\put(5544,-2587){\makebox(0,0)[lb]{\smash{\SetFigFont{6}{7.2}{\familydefault}{\mddefault}{\updefault}$5/3$}}}
\put(5094,-4036){\makebox(0,0)[lb]{\smash{\SetFigFont{6}{7.2}{\familydefault}{\mddefault}{\updefault}$10/3$}}}
\put(780,-4029){\makebox(0,0)[lb]{\smash{\SetFigFont{6}{7.2}{\familydefault}{\mddefault}{\updefault}$4/3$}}}
\put(369,-2580){\makebox(0,0)[lb]{\smash{\SetFigFont{6}{7.2}{\familydefault}{\mddefault}{\updefault}$8/3$}}}
\end{picture}
\caption{Part of $G(\Phi,1)$}
\label{k3ext}
\end{center}
\end{figure}

The polynomial of degree 2 such that $L_1(1)=2$, $L_1(2)=3$ and
$L_1(3)=1$ is $L_1(x)=-\frac{3}{2}x^2+\frac{11}{2}x-2$, and the
polynomial of degree 2 such that $L_2(1)=3$, $L_2(3)=2, L_2(2)=1$ is
$L_2(x)=\frac{3}{2}x^2-\frac{13}{2}x+8$.  If
$\Phi(x,y)=(y-L_1(x))(y-L_2(x))$, then $\vec{G}(\Phi, 1)$ is $D=K_3$.
Nevertheless, the component $G(\Phi,1)$ is infinite.  In the
Figure~\ref{k3ext} vertices adjacent to 1, 2 and 3 which are in
$G(\Phi,1)$ but not in
$\vec{G}(\Phi,1)$ are shown.

By Theorem~\ref{strong} the condition on a finite strongly connected
$d$-regular digraph of being isomorphic to a strong component of $G(\Phi)$ for
some polynomial $\Phi(x,y)$ is not restrictive at all. A $d$-regular digraph
$D$ is said to be \emph{polynomial} if, for some standard polynomial
$\Phi(x,y)$, the digraph $D$ is isomorphic to $G(\Phi)$ or to a non singular
component of $G(\Phi)$.

\section{Cayley digraphs}

Cayley digraphs are relevant structures in different contexts such
as modeling interconnection networks~\cite{Heydemann,LaJwDh},
tessellations of the sphere and of the Euclidean
Plane~\cite{Bollobas} and in combinatorial group
theory~\cite{White}. Recall that, given a group $\Gamma$ and a finite set
$S\subseteq\Gamma$ with $1\not\in S$, the \emph{Cayley digraph}
$\textrm{Cay}(\Gamma,S)$ is defined by taking the elements in
$\Gamma$ as vertices and an ordered pair $(u,v)$ is an arc if
$v=su$ for some $s\in S$. A Cayley digraph
$\textrm{Cay}(\Gamma,S)$ is connected if and only if $S$
is a generating set of $\Gamma$ (this is the reason why the
condition of $S$ being a generating system is often included in
the definition). If $s^{-1}\in S$ for all $s\in S$, then
$\textrm{Cay}(\Gamma,S)$ is a graph.  Cayley digraphs are known to
be vertex transitive. This implies that all components of a Cayley digraph are
isomorphic and also that all strong components are isomorphic.

\begin{thm}
\label{cayley+}
\begin{enumerate}
\item [\rm (i)] Let $\Phi(x,y)$ be a standard polynomial. If $\Phi(x,y)=f(y-x)$
  for some univariate polynomial $f(s)\in\mathbb{C}[s]$, then
  $G(\Phi)$ is a Cayley digraph on $(\mathbb{C},+)$.
\item [\rm (ii)] Let $D=\mathrm{Cay}(\mathbb{C},S)$ be a Cayley digraph on
  $(\mathbb{C},+)$.  Then there exists a univariate polynomial
  $f(s)\in\mathbb{C}[s]$ such that $D$ is isomorphic to $G(\Phi)$
  where $\Phi(x,y)=f(y-x)$ is a standard polynomial.
\end{enumerate}
\end{thm}
\begin{pf}
  (i) Let $s_1,\ldots,s_d$ be the roots of $f(s)$. Then
  $\Phi(x,y)=f(y-x)=c(y-x-s_1)\cdots (y-x-s_d)$ with $c\ne 0$.  In
  $G(\Phi)$ a vertex $u$ is adjacent to the vertices
  $u+s_1,\ldots,u+s_d$. Because $\Phi(x,y)$ is standard, $s_i\ne 0$
  for all $i$ and $s_i\ne s_j$ for $i\ne j$. If $S=\{s_1,\dots,s_d\}$,
  then we get $G(\Phi)=\mathrm{Cay}(\mathbb{C},S)$.

(ii) Given $\mathrm{Cay}(\mathbb{C},S)$ with $S=\{s_1,\ldots, s_d\}$,
consider the polynomial $f(s)=(s-s_1)\cdots(s-s_d)$ and take
$\Phi(x,y)=f(y-x)$. The polynomial $\Phi(x,y)$ is standard and
$\mathrm{Cay}(\mathbb{C},S)=G(\Phi)$. \qed
\end{pf}

Note that for a polynomial $\Phi(x,y)=f(y-x)$ as in
Theorem~\ref{cayley+} the components are always infinite.  For
instance, if $\Phi(x,y)=(y-x)^4-1=(y-x-1)(y-x+1)(y-x-i)(y-x+i)$ then
$G(\Phi)$ has no singular vertices and it is isomorphic to
$\textrm{Cay}(\mathbb{C},\{1,-1,i,-i\})$. In this example, components and
strong components coincide and all of them are isomorphic to
$G(\Phi,0)$, the grid of integer coordinates.

Next theorem is the corresponding to Theorem~\ref{cayley+}
for Cayley digraphs on the
multiplicative group of $\mathbb{C}$. As usual,
$\mathbb{C}^*$ denotes $\mathbb{C}\setminus\{0\}$.

\begin{thm}
\label{cayley*}
\begin{enumerate}
\item[\rm (i)] Let $\Phi(x,y)$ be an homogeneous standard polynomial.
  Then $G(\Phi)^*$ is a Cayley digraph on $(\mathbb{C}^*,\cdot)$.
\item[\rm (ii)] Let $\mathrm{Cay}(\mathbb{C}^*, S)$ be a Cayley digraph on
  $(\mathbb{C}^*,\cdot)$. Then, there exists an homogeneous standard
  polynomial $\Phi(x,y)$ such that $\mathrm{Cay}(\mathbb{C}^*,
  S)=G(\Phi)^*$.
\end{enumerate}
\end{thm}
\begin{pf}
  (i) Let $\Phi(x,y)$ be an homogeneous standard polynomial of total
  degree $d$. Note that $\Phi(x,y)$ being standard, it must be also of
  partial degree $d$ in both indeterminates. We have
  $\Phi(x,sx)=x^df(s)$ where $f(s)$ is a univariate polynomial in $s$
  of degree $d$.  Let $s_1,\ldots,s_d$ be the roots of $f(s)$. Then,
  $\Phi(x,s_ix)=0$ for $1\le i\le d$ and
  $\Phi(x,y)=c(y-s_1x)\cdots(y-s_dx)$ for some $c\ne 0$. As
  $\Phi(x,y)$ is standard, $s_i\ne 1$ and $s_i\ne 0$ for all $i$, and
  $s_i\ne s_j$ for $i\ne j$.  Each vertex $u\in\mathbb{C}^*$ is
  adjacent to the $d$ vertices $s_1u,\ldots,s_du$.  Therefore,
  if $S=\{s_1,\ldots,s_d\}$, we have
  $G(\Phi)^*=\mathrm{Cay}(\mathbb{C}^*, S)$.

  (ii) Given a Cayley digraph $\mathrm{Cay}(\mathbb{C}^*, S)$ on the
  multiplicative group on $(\mathbb{C}^*,\cdot)$, where
  $S=\{s_1,\ldots,s_d\}$, then $\Phi(x,y)=(y-s_1x)\cdots(y-s_dx)$ is a
  standard polynomial and $G(\Phi)^*=\mathrm{Cay}(\mathbb{C}^*, S)$. \qed
\end{pf}

Note that if $\Phi(x,y)$ is an homogeneous standard polynomial of total degree
$d$, then $G(\Phi)$ has 0 as the unique singular vertex, and $(0,0)$ is
a loop of multiplicity $d$.

As a Corollary of Theorems~\ref{cayley+} and~\ref{cayley*} we have:

\begin{cor}
Cayley digraphs on the additive and multiplicative groups of $\mathbb{C}$
are polynomial.
\end{cor}

A \emph{circulant} digraph is a strongly connected Cayley digraph on a
finite cyclic group.  It is not a restriction to take the group $U_n$
of the $n$-th roots of the unity as the cyclic group of order $n$.
Then, a circulant digraph is a Cayley digraph of the form
$\textrm{Cay}(U_n,S)$, where $S$ is a generating set of $U_n$.  Now,
$\textrm{Cay}(U_n,S)$ is the component of $1$ in
$\textrm{Cay}(\mathbb{C}^*,S)$, so we conclude
\begin{cor}
Circulant digraphs are polynomial.
\end{cor}
For instance, if $\omega$ is a primitive $n$-root of unity and we
define $\Phi(x,y)=\prod_{i=1}^{n-1}(y-\omega^ix)$, the components of
$G(\Phi)^*$ are complete graphs $K_{n}$. If $\omega$ is a primitive
$2d$-root of unity and $\Phi(x,y)=\prod_{i=1}^{d}(y-\omega^{2i-1}x)$,
then the components of $G(\Phi)^*$ are complete bipartite graphs $K_{d,d}$.

The $n$-prisms are a family of Cayley digraphs over non cyclic groups that
are also polynomial. The $n$-\emph{prism} is the
Cayley graph
$$
\mathrm{Cay}(\mathbb{Z}_n\times\mathbb{Z}_2,\{(1,0),(n-1,0),(0,1)\}).
$$
For instance, the $3$-dimensional cube is the $4$-prism. The $n$-prism
can be obtained by the polynomial $\Phi(x,y)=(y-\omega
x)(y-\omega^{n-1}x)(xy-2)$, where $\omega$ is a $n$-th primitive root
of the unity.

\section{Polynomials of partial degree one}

If $\Phi(x,y)=a_1(x)y+a_0(x)$ is a polynomial of partial degree one in $y$,
then the vertices $v$ such that there exists a directed path from $u$ to $v$
are the vertices $u_n$ defined by $u_0=u$ and $u_{n+1}=-a_0(u_n)/a_1(u_n)$ for
all $n\ge 0$. Thus, the structure of $G(\Phi)$ is closely related to the
dynamical system defined by the rational function $f(x)=-a_0(x)/a_1(x)$.  The
iteration of rational functions was studied by G. Julia~\cite{Julia} as early
as 1918 and P. Fatou~\cite{Fatou} in 1922, and there is a huge literature on
the topic (see, for instance \cite{Keller}). In this section we give a
complete description of the digraphs $G(\Phi)$ when $\Phi(x,y)$ is a
polynomial of partial degree one in both indeterminates.

Note that a component of a
digraph of indegree and outdegree equal to one is isomorphic either to a
directed $n$-cycle $\vec{C}_n=\mathrm{Cay}(\mathbb{Z}_n,\{1\})$ for
some $n$ or to an infinite path
$\vec{P}=\mathrm{Cay}(\mathbb{Z},\{1\})$.  We shall see that in any
case, all components of $G(\Phi)^*$ are isomorphic. This result is
applied to exhibit examples of Cayley digraphs on non commutative
groups that are polynomial digraphs.  First, we characterize the
standard polynomials of partial degree 1.
\begin{lem}
\label{standard}
  A polynomial $\Phi(x,y)=(cx+d)y-(ax+b)$ is standard if and only if
   $ad-bc\ne 0$ and it is not divisible by $y-x$.
\end{lem}
\begin{pf}

  Assume that $\Phi(x,y)=(cx+d)y-(ax+b)=(cy-a)x+dy-b$ is standard.
  Then the polynomials $D(x)=cx+d$ and $E(x)=cy-a$ are not the zero
  polynomials.  First, consider the case $c=0$. Then $d\ne 0$ and
  $a\ne 0$, so $ad-bc=ad\ne 0$. Second, assume $a=0$. Then $c\ne 0$.
  If $b=0$, then $B(y)=\mathrm{gcd}(cy-a,dy-b)=\mathrm{gcd}(cy,dy)\ne
  1$, a contradiction.  Thus, $b\ne 0$ and $ad-bc=-bc\ne 0$. Finally,
  let $ac\ne 0$.  By dividing $cx+d$ by $ax+b$, the remainder is
  $-bc/a+d=(ad-bc)/a$.  As $A(x)=\mathrm{gcd}(cx+d,ax+b)=1$, this
  remainder must be not zero.  Therefore $ad-bc\ne 0$. The condition of
  not being divisible by $y-x$ ensures that $L(x)$ is not the zero
  polynomial.

  Conversely, assume that $ad-bc\ne 0$. Then $c$ and $d$ can not be
  simultaneously zero, so $D(x)=cx+d$ is not the zero polynomial.
  Analogously, $E(x)=cy-a$ is not the zero polynomial. If
  $A(x)=\mathrm{gcd}(cx+d,ax+b)$ is non constant, then $c=a\lambda$ and
  $d=b\lambda$, which implies $ad-bc=0$, a contradiction. Analogously,
  $B(x)=\mathrm{gcd}(cy-a,dy-b)$ non constant implies $ad-bc=0$.
  Finally, if $L(x)=cx^2+(d-a)x-b$ is the zero polynomial, then $b=c=0$
  and $d=a$. Therefore $\Phi(x,y)=d(y-x)$ is divisible by $y-x$. \qed
\end{pf}

For standard polynomials $\Phi(x,y)=(cx+d)y-(ax+b)$ of partial degree
one, the structure of $G(\Phi)^*$ is --as observed above-- closely
related to the properties of the maps $f(z)=(az+b)/(cz+d)$ with
$ad-bc\ne 0$ (or, equivalently, with $ad-bc=1$). These maps, called
\emph{linear fractional transformations}~\cite{Ford} or \emph{M\"oebius
  transformations}~\cite{Caratheodory}, apply the complex plane minus
$-d/c$ to the complex plane minus $a/c$. They are examples of
conformal maps and form a group $M(\mathbb{C})$ under composition.  It
is useful to represent the group $M(\mathbb{C})$ as a quotient of the
group $S_2(\mathbb{C})$ of square matrices of order $2$ with
determinant 1 as follows. The mapping
$$
\begin{array}{rcl}
S_2(\mathbb{C}) & \longrightarrow & M(\mathbb{C}) \\
\left(\begin{array}{cc} a & b \\ c & d \end{array}\right) & \mapsto &
\displaystyle f(z)=\frac{az+b}{cz+d}
\end{array}
$$
is an surjective group homomorphism and its kernel is $\{1,-1\}$.
Therefore $ M(\mathbb{C})\simeq S_2(\mathbb{C})/\{+1,-1\}$.  Given
$f\in M(\mathbb{C})$ defined by $f(z)=(az+b)/(cz+d)$ we denote
$A_f=\left(\begin{array}{cc} a & b \\ c & d
\end{array}\right)$ and $[A_f]$ the class of $A_f$ in $S_2(\mathbb{C})/\{+1,-1\}$.

\begin{thm}
\label{bilinear}
If $\Phi(x,y)=(cx+d)y-(ax+b)$ is a standard polynomial, then all
components of $G(\Phi)^*$ are isomorphic.
\end{thm}
\begin{pf}
  Let $\ell_1, \ell_2$ be the roots of $L(x)=\Phi(x,x)=cx^2+(d-a)x-b$.
  Note that if $c=0$ all vertices have outdegree $1$; otherwise,
  $-d/c$ is the unique out-defective vertex.  Let $V=\mathbb{C}$ if
  $c=0$ and $V=\mathbb{C}\setminus\{-d/c\}$ if $c\ne 0$.  A vertex $u$
  in $V$ is adjacent to the vertex $f(u)=(au+b)/(cu+d)$.  Consider the
  function $f(x)=(ax+b)/(cx+d)$ defined on $V$.  Assume that there
  exist components of $G(\Phi)^*$ which are directed cycles and let
  $n\ge 2$ be the minimum of the lengths of these cycles. Then $n$ is
  the minimum positive integer such that there exists a vertex $u$ in a
  non-singular component such that $f^n(u)=u$, or equivalently,
 $$
u=f^n(u)=\frac{a_nu+b_n}{c_nu+d_n} \mbox{ where }
\left(\begin{array}{cc} a_n & b_n \\ c_n & d_n \end{array}\right)
=A_f^n.
$$
Thus, the solutions of $f^n(x)=x$ are the roots of
$F(x)=c_nx^2+(d_n-a_n)x-b_n$. But $\ell_1$, $\ell_2$ and $u$ are three
roots of $F(x)$ (if $\ell_1=\ell_2$, then $\ell_1$ has multiplicity
2), so $F(x)$ is the zero polynomial and $f^n(v)=v$ for all $v$. This
implies that all components of $G(\Phi)^*$ are isomorphic to
$\vec{C}_n$. \qed
\end{pf}
\begin{table}[t]%
\caption{\label{abcd} Conditions for the components of $G(\Phi)^*$ being directed $n$-cycles.}
\begin{tabular}{rl}
$n$ & $\vec{C}_n(a,b,c,d)$ \\
\hline
2 & $a + d$ \\
3 & $a^{2} + b\,c + a\,d + d^{2}$ \\
4 & $a^{2} + 2\,b\,c + d^{2}$ \\
5 & $a^{4} + 3\,a^{2}\,b\,c + b^{2}\,c^{2} + a^{3}\,d
          + 4\,a\,b\,b\,d + a^{2}\,d^{2} + 3\,b\,c\,d^{2} + a\,d^{3} + d^{4}$\\
6 & $3\,b\,c + a^{2} - a\,d + d^{2}$\\
7 & $8\,c\,a^{3}\,b\,d + 9\,c\,a^{2}\,b\,d^{2} + 6\,c^{2}\,a^{2}\,b^{2}
                       + 9\,c^{2}\,a\,b^{2}\,d + a^{6} + 8\,c\,a\,b\,d^{3}
                       + 5\,a^{4}\,b\,c + c^{3}\,b^{3}$ \\
  & $\phantom{8\,c\,a^{3}\,b\,d} + 6\,c^{2}\,b^{2}\,d^{2}
                                 + 5\,d^{4}\,c\,b + d\,a^{5}
                                 + d^{2}\,a^{4} + d^{3}\,a^{3} + d^{4}\,a^{2}
                                 + d^{5}\,a + d^{6}$\\
8 & $2\,c^{2}\,b^{2} + a^{4} + d^{4} + 4\,a^{2}\,b\,c + 4\,d\,a\,b\,c
                    + 4\,c\,b\,d^{2}$\\
9 & $c^{3}\,b^{3} + 9\,c^{2}\,b^{2}\,d^{2} + 15\,c^{2}\,a\,b^{2}\,d
                 + 9\,c^{2}\,a^{2}\,b^{2} + 6\,c\,a^{3}\,b\,d
                 + 6\,d^{4}\,c\,b + 3\,c\,a^{2}\,b\,d^{2}$\\
  & $\phantom{c^{3}\,b^{3}} + 6\,a^{4}\,b\,c + 6\,c\,a\,b\,d^{3}
                             + d^{6} + a^{6} + d^{3}\,a^{3}$ \\
10 & $5\,c^{2}\,b^{2} - d\,a^{3} + d^{2}\,a^{2} + 5\,c\,b\,d^{2}
                     + 5\,a^{2}\,b\,c - d^{3}\,a + d^{4} + a^{4}$\\
\hline\\
\end{tabular}
\end{table}
From a computational point of view, conditions on $a,b,c$ and $d$ for
the components of $G(\Phi)^*$ being directed $n$-cycles are easily
obtained. The polynomial $F(x)$ in the above proof is of the form
$F(x)=c_nx^2+(d_n-a_n)x-b_n=F_n(a,b,c,d)L(x)=F_n(a,b,c,d)(cx^2+(d-a)x-b)$.
Therefore, $F_n(a,b,c,d)=c_n/c=b_n/b=(d_n-a_n)/(d-a)$.
  Now, for each
divisor $k$ of $n$, the polynomial $F_k(a,b,c,d)$ must be a factor of
$F_n(a,b,c,d)$. By dividing $F_n(a,b,c,d)$ by all the factors
corresponding to digraphs $G(\Phi)^*$ with directed $k$-cycles as components
($k$ divisor of $n$), we obtain the condition $\vec{C}_n(a,b,c,d)=0$ for
the components of $G(\Phi)^*$ to be directed cycles of length
$n$. In the table~\ref{abcd} the conditions $\vec{C}_n(a,b,c,d)=0$ for $n$
from 2 to 10 are given.

The above proposition can be interpreted in the sense that if a linear
fractional transformation $f(z)=(az+b)/(cz+d)$ generates a cyclic
group of order $n$, then all its orbits but the fixed points (the
loops in the graph) have $n$ points (and therefore the components of
$G(\Phi)^*$, where $\Phi(x,y)=(cx+d)y-(ax+b)$, are directed $n$-cycles).

The following theorem is a characterization of polynomials of partial
degree one which give directed $n$-cycles. By a \emph{primitive}
$n$-root of a complex number $z$ we mean a complex $r$ such that $n$
is the minimum positive integer with $r^n=z$.

\begin{thm}
\label{ncycles}
  Let $\Phi(x,y)=(cx+d)y-(ax+b)$ be a standard polynomial and let $r$
  be a square root of $(a+d)^2-4$.  Then the components of $G(\Phi)^*$
  are directed $n$-cycles if and only if $a+d\not\in\{+2,-2\}$ and
  $(a+d\pm r)/2$ are both $n$-th roots of \,$1$ or both $n$-roots of $-1$,
  and at least one of them primitive.
\end{thm}
\begin{pf}
  Let $f\in M(\mathbb{C})$ defined by $f(z)=(az+b)/(cz+d)$ with
  $ad-bc=1$.
  We have the following equivalences:\\
  The components of $G(\Phi)^*$ are isomorphic to $\vec{C}_n$ \\
  $\Leftrightarrow$ $f^n$ is the identity mapping \\
  $\Leftrightarrow$ the order of $[A_f]$ is  $n$ \\
  $\Leftrightarrow$ $A_f$ has two distinct eigenvalues which are
  $n$-th roots of $1$ or $-1$ and at least one of them primitive.  The
  characteristic polynomial of $A_f$ is $\lambda^2-(a+d)\lambda+1$ and
  the discriminant $(a+d)^2-4$. Therefore $A_f$ has two distinct
  eigenvalues if and only if $a+d\ne \pm 2$. In this case the
  eigenvalues are $(a+d\pm r)/2$ and $[A_f]$ is of order $n$ if they are
  $n$-th roots of $1$ or $-1$ and at least one primitive. \qed
\end{pf}

For instance, consider the polynomial
$\Phi(x,y)=(x-1+\sqrt{3})y-(x-2+\sqrt{3})$. We have the values
$a=c=1$, $b=-2+\sqrt{3}$ and $d=-1+\sqrt{3}$. They satisfy the 6
cycle-condition $3bc+a^2-ad+d^2=0$ of table~\ref{abcd}. On the other
hand, we can use the Theorem~\ref{ncycles}. The matrix associated to the
corresponding $f$ is $A_f=\left(\begin{array}{cc}
    1\phantom{x} & -2+\sqrt{3} \\
    1\phantom{x}  & -1+\sqrt{3}
\end{array}\right)$. The eigenvalues of $A_f$ are $\lambda=(+\sqrt{3}\pm i)/2$
which are primitive $6$-roots of $-1$. Therefore $[A_f]$ has order 6 and
all components of $G(\Phi)^*$ are directed 6-cycles.

Taking into account that $M(\mathbb{C})$ is a non commutative
group, we can find polynomial Cayley digraphs on non abelian
groups. The following theorem opens the way.

\begin{thm} Let $\Gamma$ be a subgroup of $M(\mathbb{C})$ generated
by a set $S=\{f_1,\ldots,f_d\}$ of $d$ linear fractional transformations.
Then $\mathrm{Cay}(\Gamma,S)$ is polynomial.
\end{thm}
\begin{pf}
  Let $f\in\Gamma$ be a transformation defined by
  $f(z)=(az+b)/(cz+d)$. We associate to $f$ the polynomial
  $\Phi_f(x,y)=(cx+d)y-(ax+b)$. The group $\Gamma$ is generated by a
  finite set, so it is a countable set. Therefore, the set of
  $u\in\mathbb{C}$ such that there exists $f\in \Gamma$ with $u$
  belonging to a singular component of $G(\Phi_f)$ is also a countable
  set. Thus, we can choose $u\in\mathbb{C}$ such that
  $\vec{G}(\Phi_f,u)$ is not a singular component of $G(\Phi_f)$ for
  all $f\in\Gamma$.

  Let $\Phi(x,y)=\Phi_{f_1}(x,y)\cdots\Phi_{f_d}(x,y)$. The mapping
  $\mathrm{Cay}(\Gamma,S)\rightarrow \vec{G}(\Phi,u)$ defined by $f\mapsto
  f(u)$ is a digraph isomorphism.  Indeed: it is injective, because $f(u)=u$
  implies that $f$ is the identity (otherwise $u$ would be a loop vertex in
  $G(\Phi_f)$). It is surjective, because if $v$ is a vertex in
  $\vec{G}(\Phi,u)$, then there exists a path $u=u_0,\ldots, u_\ell=v$.  Now,
  each arc is of the form $(u_j,f_{i_j}(u_j))$ for some $f_{i_j}\in S$. If
  $f=f_{i_\ell}\cdots f_{i_1}$, we have $f(u)=v$.  Finally, it preserves
  adjacencies: $(g,h)$ is an arc in $\mathrm{Cay}(\Gamma, S)$
  $\Leftrightarrow$ $g=f_ih$ for some $f_i\in S$ $\Leftrightarrow$
  $g(u)=f_i(h(u))$ for some $i$ $\Leftrightarrow$ $\Phi_{f_i}(h(u),g(u))=0$
  for some $i$ $\Leftrightarrow$ $\Phi(x,y)=0$. \qed
\end{pf}
The finite subgroups of $M(\mathbb{C})$ are determined (\cite{Ford},
Chapter VI). In particular, dihedral groups and the groups of
symmetries of regular polyhedra are finite subgroups of $M(\mathbb{C})$.
Therefore, we have:
\begin{cor}
  Cayley digraphs on dihedral groups and on the groups of symmetries of
  regular polyhedra are polynomial.
\end{cor}
For instance, the dihedral group $D_{2n}=\langle f,t\  |\
f^n=t^2=1, \ tftf=1\rangle$ is the subgroup of $M(\mathbb{C})$ generated
by $f(z)=\omega z$ and $t(z)=2/z$, where $\omega$ is a primitive
$n$-root of unity. Then, the Cayley digraph
$\mathrm{Cay}(D_{2n},\{f,t\})$ is obtained by the polynomial
$\Phi(x,y)=(y-\omega x)(xy-2)$.

\section{Symmetric polynomials of degree two}
In this section we give a method for analyzing the components of $G(\Phi)$ for
a standard symmetric polynomial $\Phi(x,y)$ of total degree two. The method
implies long but routine calculations, so we skip them and give only
results.

From
Theorem~\ref{bilinear} we can assume that the standard polynomial is not of
partial degree one, so it is of the form $\Phi(x,y)=x^2+y^2+axy+b(x+y)+c$.
Since $G(\Phi)^*$ is a 2-regular graph, a component of $G(\Phi)^*$ is
isomorphic to a (undirected) cycle $C_n=\mathrm{Cay}(\mathbb{Z}_n,\{1,-1\})$
or to the double ray graph $R=\mathrm{Cay}(\mathbb{Z},\{1,-1\})$.
Denote $\ldots,v_{-n},\ldots,v_{-1},v_0,v_1,\ldots,v_n,\ldots$ the vertices of
$G(\Phi,v_0)$ with $v_i$ and $v_{i+1}$ adjacent vertices.
The polynomial of degree two $\Phi(v_{n-1},y)$
has two roots, namely $v_{n-2}$ and $v_n$.
As the sum of the two roots is minus de coefficient of $y$, we have
$v_{n-2}+v_n=-(av_{n-1}+b)$ or, equivalently,
\begin{equation}
\label{recurrence}
v_n+av_{n-1}+v_{n-2}=-b.
\end{equation}
The two roots of $\Phi(v_0,y)$ are $v_{-1}$ and $v_1$.
The solutions of the recurrence~(\ref{recurrence}) with initial values $v_0$
and $v_1$ determine the vertices with positive subscripts and analogously for
the negative subscripts taking $v_0$ and $v_{-1}$ as initial values.
Then, to determine
the structure of $G(\Phi)$ the method is the following. First, to solve the
recurrence~(\ref{recurrence}), wich give $v_n$ in terms of two initial values
$v_0$ and $v_1$. Second, to find the singular vertices. There are not
defective vertices, so only vertices with loops and multiple arcs should be
calculated. In any case, there exist at most two loops
and two origin of multiple arcs. For each loop-vertex $\ell$, the solution of
(\ref{recurrence}) for the initial values $v_0=v_1=\ell$ gives the
vertices in $G(\Phi,\ell)$. Analogously, for each multiple arc $(m,m_1)$,
the solution of~(\ref{recurrence}) for the initial conditions $v_0=m$ and
$v_1=m_1$ gives the vertices in $G(\Phi,m)$. The explicit form of the
solutions allows to decide if some of these singular components coincide and
if they are finite or not. Finally, for $v_0$ in a non singular component,
the solutions of~(\ref{recurrence})
for the initial values $v_0$ and $v_1$ and $v_0$ and $v_{-1}$ gives
the vertices in the non singular component $G(\Phi,v_0)$.

The characteristic equation of the second order linear
recurrence~(\ref{recurrence}) is $\lambda^2+a\lambda+1=0$, and the
discriminant is $\Delta(a)=a^2-4$. To solve the recurrence three cases have to
be considered: $a=-2$, $a=2$, and $a^2-4\ne 0$. We sumarize the discussion in
each case.

\begin{figure}[t]
\caption{\label{menosdos} Structure of $G(\Phi)$ for $a=-2$}
\begin{center}
\epsfig{file=fig3.eps}
\end{center}
\end{figure}

First consider the case $a=-2$, see Figure~\ref{menosdos}.
We have $\Phi(x,y)=(x-y)^2+b(x+y)+c$.  As
$\Phi(x,y)$ is standard, $b$ and $c$ cannot be simultaneously zero.
For $b=0$ there
are no singular components. For $b\ne 0$ there exist two singular components,
both infinite, one containing the vertex-loop $\ell=-c/(2b)$ and one
containing the origin of a double arc $m=(b^2-4c)/(8b)$.
In both cases ($b=0$ and $b\ne 0$) all non singular
components are isomorphic to $R$.

For $a\ne -2$, it follows from Lemma~\ref{iso} that if
$\Psi(x,y)=\Phi(x-b/(a+2), y-b/(a+2))=x^2+y^2+axy+(c-b^2/(a+2))$, then
$G(\Phi)$ is isomorphic to $G(\Psi)$. Therefore we can assume (and we
do) without loss of generality that the polynomial $\Phi(x,y)$ is of the form
$\Phi(x,y)=x^2+y^2+axy+c$.

\begin{figure}[h]
\caption{\label{dos} Structure of $G(\Phi)$ for $a=2$}
\begin{center}
\epsfig{file=fig4.eps}
\end{center}
\end{figure}

Consider now the case $a=2$, see Figure~\ref{dos}.
Now the polynomial is $\Phi(x,y)=(x+y)^2+c$.  As
$\Phi(x,y)$ is a standard polynomial, necessarily $c\ne 0$. There exist two
non-singular components, both infinite, one containing a loop at the vertex
$\ell_1=\sqrt{-c}/2$, where $\sqrt{-c}$ is one of the two square roots of
$-c$, and the other containing a loop at the vertex $\ell_2=-\ell_1$.
The non singular components are isomorphic to $R$.

Finally, consider the case $a^2-4\ne 0$, see Figures~\ref{noroot}
and~\ref{root}.
If $c=0$, there exists only one singular component, which contains only the
vertex $\ell=0$ with a double loop. If $c\ne 0$, there exist two
loop-vertices
$\ell_1=\sqrt{-c/(a+2)}$ and $\ell_2=-\ell_1$ and two origin of double arcs,
$m_1=2\sqrt{c/(a^2-4)}$ and $m_2=-m_1$. The number and the finitness of the
components depends on the values of $a$ or, equivalently, on the values of
the two distinct roots $\omega_1$ and
$\omega_2=1/\omega_1$ of $\lambda^2+a\lambda+1$.
If $\omega_i$ are primitive $n$-th roots of the unity for some positive
integer $n$, then the components of $G(\Phi)$ are finite, otherwise they are
infinite. Now, $\omega_i$ is a primitive $n$-th root of the unity if and only
if $a$ is of the form $a=2\cos(2\pi k/n)$ for some $k$ with $\gcd(k,n)=1$.

\begin{figure}[t]
\caption{\label{noroot} Structure of $G(\Phi)$ for $a^2-4\ne 0$, $a\ne
  2\cos\frac{2\pi k}{n}$, $\gcd(k,n)=1$}
\begin{center}
\epsfig{file=fig5.eps}
\end{center}
\end{figure}

The following theorem sumarizes the form of the non singular components.

\begin{thm}
\label{thm degree 2}
Let $\Phi(x,y)=x^2+y^2+axy+b(x+y)+c$ be a standard polynomial. If $a\ne\pm 2$
and $a=2\cos(2\pi k/n)$ for some positive integers $n$ and $k$ with
$\gcd(k,n)=1$, then all components of $G(\Phi)^*$ are isomorphic to $C_n$, the
cycle of length $n$.  Otherwise, they are isomorphic to the double ray graph
$R$.
\end{thm}

\begin{figure}[t]
\caption{\label{root} Structure of $G(\Phi)$ for $a^2-4\ne 0$,
$a=2\cos\frac{2\pi k}{n}$, $\gcd(k,n)=1$}
\begin{center}
\epsfig{file=fig6.eps}
\end{center}
\end{figure}

\section{A conjecture}

If $\Phi(x,y)$ is an homogeneous standard polynomial, then $G(\Phi)^*$
is a Cayley digraph (Theorem~\ref{cayley*}). Therefore if a component
of $G(\Phi)^*$ is finite, all of them are finite and isomorphic. In
Section~5 we have seen that if $\Phi(x,y)$ is a polynomial of partial
degree one, then if a component of $G(\Phi)^*$ is a directed $n$-cycle,
then all of them are directed $n$-cycles.  Also, for a symmetric
polynomials $\Phi(x,y)$ of partial and total degree 2, if a component
of $G(\Phi)^*$ is a $n$-cycle, then all of them are $n$-cycles. All
these examples suggest the following conjecture:

\begin{conj}
\label{cd}
Let $\Phi(x,y)$ be a standard polynomial of partial degree $d$ and $H$ a
$d$-regular digraph isomorphic to a component of $G(\Phi)^*$.
Then, all components of $G(\Phi)^*$ are isomorphic to $H$.
\end{conj}

In~\cite{BrMo2} more evidence of Conjecture~\ref{cd} is given. For
instance, it is shown that if for a symmetric polynomial $\Phi(x,y)$ of
partial degree two $G(\Phi)^*$ has a component which is a $n$-cycle (for
small values of $n$) then all of them are $n$-cycles. Also, if
$\Phi(x,y)$ is a polynomial such that $G(\Phi)^*$ has a component
isomorphic to $K_n$ ($2\le n\le 6$) then all of them are isomorphic to
$K_n$.

\section{Acknowledgements}
We would like to thank Pelegr\'{\i} Viader for his many helpful comments
and his insightful perusal of our first draft, and Jos\'e Luis Ruiz, who first
studied $G(\Phi)$ for an homogeneous polynomial $\Phi(x,y)$ and whose
observations lead to Theorem~\ref{cayley*}.

\bibliographystyle{abbrv}

\end{document}